\documentclass[a4paper]{amsart}
\usepackage{graphicx}
\usepackage{amssymb, amsthm}
\usepackage{enumitem, bm}

\usepackage{graphicx}
\usepackage[hidelinks]{hyperref}

\theoremstyle{plain}

\theoremstyle{definition}

\newcommand{\on}{\operatorname}

\newcommand{\R}{\mathbb{R}}

\newcommand{\C}{\mathbb{C}}

\DeclareMathOperator*{\argmax}{arg\,max}
\DeclareMathOperator*{\argmin}{arg\,min}

\title[The geometry of the maximum likelihood]{The geometry of the maximum likelihood  of Cauchy-like distributions}
\author{Pavol Ševera}
\email{pavol.severa@gmail.com}
\begin{document}

\begin{abstract}
A simple way of obtaining robust estimates of the ``center'' (or the ``location'') and of the ``scatter'' of a dataset is to use the maximum likelihood estimate with a class of heavy-tailed distributions, regardless of the ``true'' distribution generating the data.
We observe that the maximum likelihood problem for the Cauchy distributions, which have particularly heavy tails, is geodesically convex and therefore efficiently solvable (Cauchy distributions are parametrized by the upper half plane, i.e.\ by the hyperbolic plane). Moreover, it has an appealing geometrical meaning: the datapoints, living on the boundary of the hyperbolic plane, are attracting the parameter by unit forces, and we search the point where these forces are in equilibrium.

This picture generalizes to several classes of multivariate distributions with heavy tails, including, in particular, the multivariate Cauchy distributions. The hyperbolic plane gets replaced by  symmetric spaces of noncompact type.
Geodesic convexity gives us an efficient numerical solution of the maximum likelihood problem for these distribution classes. This can then be used for robust estimates of location and scatter, thanks to the heavy tails of these distributions.
\end{abstract}

\maketitle

\tableofcontents

\newpage

\section{Introduction}

The Cauchy distribution
has particularly heavy tails (so heavy that $\on E[X]$ doesn't exist), so it's a natural tool for robust statistical models. It is given by the density
\begin{equation}\label{1dcauchy}
f(x\mid z)=\frac1\pi\,\on{Im}\frac1{x-z} = \frac1\pi \frac{v}{(x-u)^2 + v^2}\qquad (z=u+i\,v)
\end{equation}
where the parameter $z$ is from the upper-half plane
$$H^2=\{z\in\C\mid\on{Im}z >0\}$$
($u=\on{Re}z$ is the center and $v=\on{Im}z$ the width of the distribution).

Given data points $x_1,\dots, x_n \in\R$, let us try to find the maximum likelihood parameter $\hat z$
\begin{equation}\label{zhat}
\hat z = \argmax\limits_{z\in H^2} \sum_k \log f(x_k\mid z).
\end{equation}
(We do not necessarily suppose that the data is drawn from a Cauchy distribution, we just want to determine $\hat z$.)

As it turns out, this problem very nicely relates to the hyperbolic geometry of $H^2$. Namely, if we imagine ropes (hyperbolic straight lines) from $x_k$'s to $z$, each pulling $z$ with the unit force, then $z=\hat z$ iff the the total force acting on $z$ is 0.

\begin{figure}[h]
\includegraphics{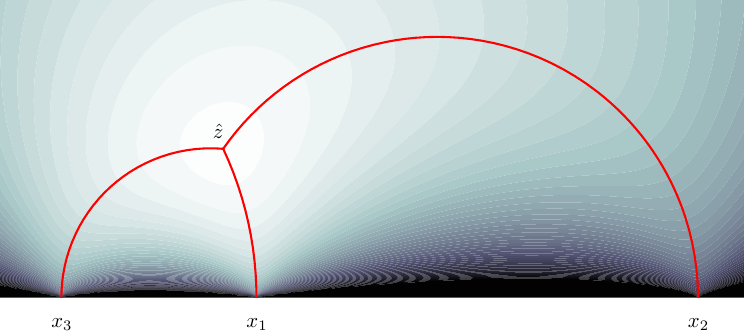}
\caption{Each datapoint $x_i$ is pulling $\hat z$ with a unit force by a red rope (hyperbolic geodesic) and these forces cancel out. The negative log likelihood (shaded) is geodesically convex.}
\end{figure}

 Moreover, the negative log likelihood $-\sum_k \log f(x_k\mid z)$ is a convex function of $z$ w.r.t.\ the hyperbolic geometry of $H^2$ (it is the sum of the (regularized) lengths of the ropes), so the solution $\hat z$ is (typically) unique and can be found by moving along hyperbolic straight lines in the direction of the resulting force (i.e.\ by geodesic gradient descent).

This picture, with geodesics connecting the data points with $\hat z$, equilibrium of unit forces and geodesic convexity of the negative log likelihood, holds also for several classes of multivariate heavy-tailed distributions, including, in particular, the multivariate Cauchy distribution. These classes are given by pairs $P\subset G$, where $G$ is a semisimple non-compact Lie group and $P$ its parabolic subgroup. The hyperbolic plane $H^2$ then gets replaced by the symmetric space $G/K$, where $K\subset G$ is a maximal compact subgroup.%
\footnote{One can follow this paper without knowing anything about these notions (except for the very last \S\ref{sec:general}); they just provide a unified picture and  explanation of our examples. On the other hand, for readers who get overwhelmed by formulas (like me), it is nice to know that no formulas are needed.}
\medskip

From a practical/statistical point of view, the problem \eqref{zhat} is simply trying to find a location $\on{Re}\hat z$ and scatter $\on{Im}\hat z$ of the dataset $x_1,\dots,x_n$ in a way that is robust with respect to  outliers (only seldom we would have reasons to believe that the data is actually drawn from a Cauchy distribution).%
\footnote{In a more proper statistical language: If $x_1,\dots, x_n$ are iid random variables with probability measure $P$ on $\R$, then $\hat z$ (given by \eqref{zhat}) is an estimator of
$$z_P:=\argmax_{z\in H^2}\on{E}_P\bigl[\log f(X\mid z) - \log f(X\mid z_0)\bigl]$$ 
where $z_0\in H^2$
is arbitrarily chosen (the subtracted term is $z$-independent and it just ensures that the expectation exists for any $P$). Assume that $P$ doesn't have an atom of probability $\geq1/2$. Then the function $z\mapsto \on{E}_P\bigl[\log f(X\mid z) - \log f(X\mid z_0)\bigl]$ is strictly geodesically convex and goes to $\infty$ at the ideal boundary of $H^2$. This ensures, in particular, the existence of $z_P$ for any such $P$. Using the methods of \cite{Huber} it then follows that the estimator $\hat z$ converges to $z_P$ in probability for any such  $P$.

For $P=N(\mu,\sigma^2)$ we have $z_P=\mu + i\,\lambda\,\sigma$ with $\lambda=0.612\dots$. For any $P$ symmetric w.r.t.\ a center $\mu\in\R$ we have $\on{Re}z_P=\mu$.\nopagebreak
}
 In this sense it is analogous to the more popular use of the Huber loss
$$\on{Hub}_\sigma(x)=
\begin{cases}
x^2/2 & (|x|\leq\sigma)\\
\sigma\,x - \sigma^2/2 & (|x|>\sigma)
\end{cases},
$$
which is a convex  and therefore convenient function. We simply point out that the Cauchy distribution, which gives the negative log likelihood loss
$$
\log(\sigma+x^2/\sigma)
$$
(which grows only logarithmically as $x\to\infty$ and so is even less sensitive to outliers),
 and also  its various multivariate relatives, lead to (geodesically) convex problems too, with a nice geometrical meaning.
 
From the multivariate distributions that we consider we also get an estimate of the ``location'' and of the ``scatter'' of the data set, but the meaning of ``scatter'' depends on the type of the distribution. In particular, in the multivariate Cauchy case it is an arbitrary ellipsoid (and so it can be used to estimate the covariance matrix of e.g.\ normal distributions contaminated by outliers), rather than just a simple scale. On the other hand, another distribution type, which we call ``conformal'', gives only a scale. One should make a choice which is appropriate for the data at hand.

It should be noted that the maximum likelihood estimate with multivariate Cauchy distributions has been used in robust statistics for  decades, with theoretical underpinning given in \cite{KT} (where the reader can find further motivation and context). We bring to it a geometric interpretation, which is not just pretty but also practical; in particular the geodesic convexity leads to an efficient algorithm.
\medskip

The paper is organized as follows. In \S\ref{sec:uni} we consider the case of (univariate) Cauchy distributions and of the hyperbolic plane, to gain intuition on this ``baby example''. In \S\ref{sec:regression} we briefly mention the corresponding regression problems. Then in \S\ref{sec:higher-dim} we discuss several multivariate generalizations, with an emphasis on multivariate Cauchy distributions \S\ref{sec:mv_cauchy}. Finally, in \S\ref{sec:general} we give the general construction given by a parabolic subgroup of a semisimple Lie group. In Appendix \ref{sec:appendix} we show numerical experiments with the algorithm coming from the case of multivariate Cauchy distributions (an implementation of the algorithm can be found \href{https://gist.github.com/mrwth/0e7a3885d970d9d87b0123afad2227d0}{\texttt{here}}).

This paper started as an informal note of 2 pages. Even though it got expanded over time, I decided to keep the style somewhat informal (instead of the usual definition-theorem-proof style) -- in the end, regardless of the utility of the method, it is about a cute geometric interpretation of some statistical problems.
\medskip

\noindent \textbf{Acknowledgements.}
I am grateful to Sylvain Sardy for a stimulating discussion and for a lot of encouragement. I am also grateful to the people at the Swiss Data Science Center, in particular to Guillaume Obozinski. While the 6 month internship that I spent there didn't really turn me to a data scientist, it was a great time and this paper stemmed directly from Guillaume's explanations about robust statistics.
I would also like to thank him for his questions and suggestions, which improved this paper significantly.
\medskip

\noindent \textbf{Important update.}
After this paper was posted on arXiv, I learned about the reference \cite{FR} where the probability measures on $G/P$ are discussed and, in particular, the maximum likelihood parameter is treated via Busemann functions, just as here. I still hope that this paper retains some value through its emphasis on concrete examples and algorithms.

\section{The case of (univariate) Cauchy distribution}\label{sec:uni}

\subsection{A quick and dirty review of the hyperbolic plane}

The hyperbolic metric in $H^2$ is given infinitesimally by
$$ds_\text{hyp} = \frac{|dz|}{\on{Im} z}.$$
The direct isometries of $H^2$ are the linear fractional transformations
$$z\mapsto\frac{az+b}{cz+d},$$
$a,b,c,d\in\R$, $ad-bc=1$.
The geodesics are the half-circles orthogonal to the real axis, and the vertical half-lines. In particular,
$$z(t) = i\exp t$$
is a geodesic parametrized by length;  to get a similar parametrization of any other geodesic, we can apply a suitable isometry to this one.

The elements of $\R\cup\{\infty\}$ are called ``points at infinity'' of the hyperbolic plane. Cauchy distributions live on $\R\cup\{\infty\}$ and their class is invariant under hyperbolic isometries (which act also on the boundary $\R\cup\{\infty\}$ of $H^2$); the parameter $z\in H^2$ is simply moved by the isometry \cite{Knight, McC}. 

If $x\in\R\cup\{\infty\}$ and $z\in H^2$, the hyperbolic distance between $x$ and $z$ is infinite. But for a fixed $x$ one can ``subtract infinity'' to get a reasonable value of $d_\text{hyp}(x,z)$, though natural only up to an $x$-dependent additive constant.%
\footnote{It is known under the name \emph{Busemann function} and it  is defined as 
$$\lim_{t\to\infty}d_\text{hyp}\bigl(\gamma(t), z\bigr)-d_\text{hyp}\bigl(\gamma(t), z_0\bigr),$$
 where $\gamma$ is any geodesic ending at $x$ and $z_0\in H^2$ is an arbitrarily chosen point; if we choose a different $z_0$, the limit changes by a $z$-independent constant.}
 We have
$$d_\text{hyp}(x,z) =%
\begin{cases}
 -\log f(x\mid z) +C & \text{ for } x\neq\infty\\
 -\on{Im}z + C & \text{ for } x=\infty
\end{cases}$$
where $C$ is an arbitrary $x$-dependent constant. These ``distances'' are invariant under hyperbolic isometries, but only up to $x$-dependent constants.

For a fixed $x$ the (hyperbolic) gradient of the function $z\mapsto d_\text{hyp}(x,z)$ is the unit vector field tangent to the geodesics coming from $x$. (The level curves of $z\mapsto d_\text{hyp}(x,z)$ are the horocycles, i.e. the circles tangent to $\R$ at $x$.)

\subsection{The maximum likelihood  point and equilibrium of unit forces}
Given elements $x_1,\dots,x_n\in\R\cup\{\infty\}$, let us try to find the point $\hat z\in H^2$ given by
$$\hat z = \argmin_{z\in H^2}\sum_k d_\text{hyp}(x_k,z).$$
Since we know the gradient of the functions $z\mapsto d_\text{hyp}(x_k,z)$, we get the following characterization: \emph{if $v_k$ is the unit vector tangent at $\hat z$ to the geodesic going from $\hat z$ to $x_k$, then $\sum_k v_k =0$}. 

To see that such a point is unique and that we can use gradient descent to find it (i.e.\ by moving along geodesics with the velocity $\sum_k v_k$), we should look at the (geodesic) convexity of the function we are minimizing.

\subsection{Convexity}
Let $\gamma(t)$ be a geodesic parametrized by length and let $x$ be a point at infinity. Then a quick calculation gives
$$d_\text{hyp}(x,\gamma(t)) = \log(c_1\,e^{t} + c_2\,e^{-t})$$
for a suitable $c_1, c_2\geq 0$, $c_1+c_2>0$. It is a convex function of $t$ and its 2nd derivative is in $[0,1]$. In other words, the function $z\mapsto d_\text{hyp}(x,z)$ on $H^2$ is geodesically convex.

Therefore, if $x_1,\dots,x_n\in\R\cup\{\infty\}$, the function
$$z\mapsto \sum_{i=1}^n d_\text{hyp}(x_i,z)$$
is geodesically convex, with the 2nd derivative  in $[0,n]$. If $n\geq3$ (and say if the $x_i$'s are distinct), this function goes to infinity at the infinite points of $H^2$ and is strictly geodesically convex, and therefore has a unique minimum $\hat z$ - a fact shown by a direct calculation in \cite{Cop}.

To find $\hat z$, we can use the gradient descent, moving along  geodesics. Since the 2nd derivative is in $[0,n]$, a safe step is $1/n$.

It should be noted that a (quite different) provably convergent procedure for the calculation of the maximul likelyhood estimate for Cauchy distributions was found already in \cite{OO}.

\subsection{Regression}\label{sec:regression}
Let $M$ be some space (typically $M=\R^m$, or its subset), let $m_1,\dots,m_k$ be a finite sequence of points in $M$, and let $x_1,\dots,x_k\in\R$ be ``observations done at $m_i$'s''. We are now looking for the function $\hat h\colon M\to H^2$ given by
\begin{equation}\label{regression}
\hat h = \argmin_{h\in\mathcal F}\biggl(\sum_i d_\text{hyp}\bigl(x_i,h(m_i)\bigr) + \on{F}(h)\biggr).
\end{equation}
Here $\mathcal F$ is some class of functions $M\to H^2$ and $\on{F}\colon\mathcal F\to\R$ is a functional. It is natural (though not necessarily useful) to demand  both $\mathcal F$ and $\on{F}$ to be invariant under hyperbolic isometries of $H^2$.

The probabilistic interpretation of \eqref{regression} is that we imagine the datapoints $x\in \R$ at $m\in M$ to be drawn from the Cauchy distribution $f\bigl(x\mid h(m)\bigr)$, where $h$ is some function $h\colon M\to H^2$, and we determine $\hat h$ using maximum likelihood, with $\mathcal F$ and $F$ coming from a prior imposed on $h$'s. But again, this interpretation should not be taken too seriously - we are simply describing a regression which is robust w.r.t\ presence of outliers.

The problem \eqref{regression} is more demanding than the one we had before, but the picture of ropes coming from $x_i$'s and pulling with unit forces stays in place.

Let us discuss just the case of ``1st-order splines''. We take $M=\R$ and  $\mathcal F$ all continuous and piecewise-$C^1$ maps such that
$$\on{F}(h) =\frac\alpha2\int_{-\infty}^\infty\Vert h'(t)\Vert_\text{hyp}^2\, dt$$
is finite.
Then the function $\hat h$ looks as follows (we shall use notation $t_1,\dots, t_k$ instead of $m_1,\dots, m_k$ and suppose that $t_1<\dots<t_k$). For $t<t_1$ and for $t>t_k$ it is constant. On each interval $[t_i, t_{i+1}]$ it is a geodesic arc. Finally, at each junction $t_i$ we  have the equilibrium of forces
$$\alpha\,\hat h'(t_i+0) - \alpha\,\hat h'(t_i-0) + v_i =0$$
where $v_i$ is the unit force pulling $\hat h(t_i)$ towards $x_i$.

\section{Higher-dimensional generalizations}\label{sec:higher-dim}

There are several generalizations of our story to multivariate distributions. The general picture with data points living on the boundary, geodesics connecting the data points with the parameter, equilibrium of unit forces, geodesic convexity of the negative log likelihood and a large symmetry group stays the same, just $H^2$ gets replaced by other Riemannian manifolds (symmetric spaces of non-compact type). We shall first describe some examples and then we explain where all these examples come from.

We describe the needed geometry in sufficient detail;  for an  in-depth discussion we can send the reader e.g.\ to \cite{Ebe}.

\subsection{A conformal class of distributions}\label{sec:conf}
This example simply replaces $H^2$ by the upper-half space $H^{n+1}\subset\R^{n+1}$ with its hyperbolic geometry, so we shall discuss it only briefly.

For $a>0$ and $\mathbf{b}\in\R^n$, let us consider the probability distribution on $\R^n$
$$f_n^\text{conf}(\mathbf{x}\mid a,\mathbf{b})\propto\bigl(1+\Vert (\mathbf{x} - \mathbf{b})/a\Vert^2\bigr)^{-n}.$$
 The parameters $(a,\mathbf{b})$ can be seen as points of $H^{n+1} = (0,\infty)\times\mathbb R^n$. 

The ideal boundary of $H^{n+1}$ is $\R^n\cup\{\infty\}$ and our probability distributions live on this boundary. The class of these distributions is closed under the group $SO(n+1,1)$ of the hyperbolic isometries of $H^{n+1}$, i.e.\ under the conformal transformations of  $\R^n\cup\{\infty\}$.

We have exactly the same interpretation and convexity of the negative log likelihood as in the case of $H^2$ and $\R\cup\{\infty\}$, and hence the same picture of equilibrium of unit forces at the maximum likelihood point of $H^n$.

\subsection{Multivariate Cauchy distributions}\label{sec:mv_cauchy}
We shall discuss this example in detail. Numerical experiments with the resulting algorithm are shown in Appendix \ref{sec:appendix}.
\subsubsection{Multivariate Cauchy distributions and their parametrization}
Let us consider the distributions on $\mathbb R^n$ of the form
$$f_n(\mathbf{x}\mid \mathbf{S},\mathbf{b})\propto\bigl(1+(\mathbf{x}-\mathbf b)^T \mathbf{S}^{-1}\, (\mathbf{x}-\mathbf b)\bigr)^{-(n+1)/2}, $$
where $\mathbf{S}$ is a symmetric positive-definite matrix and $\mathbf{b}\in\R^n$ (i.e.\ the $n$-variate $t$-distributions with 1 degree of freedom).
For our purposes it is more useful to parametrize these distributions a bit differently: Let $\mathbf{\tilde x}\in\mathbb R^{n+1}$ be obtained from $\mathbf{x}\in\mathbb R^{n}$ by adding $1$ as the last component. Then our distributions are
\begin{equation}\label{mv-cauchy-inv-form}
f_n(\mathbf{x}\mid \mathbf{T})\propto\bigl(\mathbf{\tilde x}^T\mathbf{T}\,\mathbf{\tilde x}\bigr)^{-(n+1)/2}
\end{equation}
where the parameter $\mathbf{T}$ is a symmetric positive definite  $(n+1)\times(n+1)$ matrix such that
$$\det\mathbf{T}=1.$$
 Importantly, the normalization constant in \eqref{mv-cauchy-inv-form} is \emph{independent of the parameter} $\mathbf{T}$, so it can be, for most purposes, ignored.
 
Let us extend  $\R^n$ to the projective space $\mathbb{RP}^n$ (which adds a set of measure 0 to $\R^n$: the elements of $\mathbb{RP}^n$ are non-zero vectors in $\R^{n+1}$ modulo the relation $\mathbf{v}\sim c\mathbf{v}$ $\forall c\neq0$, and $\mathbf{x}\in\R^n$ corresponds to the equivalence class $[\mathbf{\tilde x}]$ containing $\mathbf{\tilde x}$). The class of  multivariate Cauchy distributions is closed under the action of the group of projective transformations $\mathit{SL}(n+1)$ on $\mathbb{RP}^n$: the action of $\mathbf{A}\in \mathit{SL}(n+1)$ replaces the parameter $\mathbf{T}$ by $\mathbf{A}^T\mathbf{T}\,\mathbf{A}$.

\subsubsection{The space of positive-definite matrices of determinant 1}
Our parameter space (replacing the hyperbolic plane $H^2$, or the hyperbolic space $H^{n+1}$ from \S\ref{sec:conf}) is the space
$$\mathit{Pos}^{n+1}_1$$
of positive-definite symmetric $(n+1)\times(n+1)$ matrices of determinant 1. As we already mentioned, the group $\mathit{SL}(n+1)$ acts (transitively) on $\mathit{Pos}^{n+1}_1$ via 
$$\mathbf{T\mapsto A^\mathit{T}T\,A}\qquad
(\mathbf{T}\in\mathit{Pos}^{n+1}_1,\ \mathbf{A}\in \mathit{SL}(n+1)).$$

The space $\mathit{Pos}^{n+1}_1$ has a natural $\mathit{SL}(n+1)$-invariant Riemannian metric: if $\mathbf{V}$ and $\mathbf{W}$ are tangent vectors at a point $\mathbf{T}\in\mathit{Pos}^{n+1}_1$, i.e.\ symmetric matrices such that 
$$\on{Tr}\bigl(\mathbf{T}^{-1}\mathbf{V}\bigr)=\on{Tr}\bigl(\mathbf{T}^{-1}\mathbf{W}\bigr)=0,$$
 their inner product is 
$$(\mathbf{V}, \mathbf{W})=\on{Tr}\bigl(\mathbf{T}^{-1}\mathbf{VT}^{-1}\mathbf{W}\bigr).$$
The geodesic $\gamma$ such that $\gamma(0)=\mathbf{T}$, $\gamma'(0)=\mathbf{V}$ is given by the action of $\exp\bigl(\frac t2\, \mathbf{T}^{-1}\mathbf{V}\bigr)$ on $\mathbf{T}$, i.e.
\begin{equation}\label{pos-geod}
\gamma(t)=\exp\bigl(\tfrac t2\, \mathbf{T}^{-1}\mathbf{V}\bigr)^T\,\mathbf{T}\,
\exp\bigl(\tfrac t2\, \mathbf{T}^{-1}\mathbf{V}\bigr).
\end{equation}
If $\gamma(0)=\mathbf{I}$ (the identity matrix), this simplifies to $\gamma(t)=\exp(t\,\mathbf{V})$, where $V$ is now a traceless symmetric matrix.

\subsubsection{The ideal boundary of $\mathit{Pos}^{n+1}_1$}
Just as $H^n$, the space $\mathit{Pos}^{n+1}_1$ has an ideal boundary (composed of ``points at infinity''). Each geodesic ends at an ideal point, and, by definition, two geodesics $\gamma_1$, $\gamma_2$ with $\Vert\gamma_1'(0)\Vert=\Vert\gamma_2'(0)\Vert$  end at the same ideal point iff the distance $d\bigl(\gamma_1(t), \gamma_2(t)\bigr)$ stays bounded as $t\to\infty$. Geodesics $t\mapsto\gamma(t)$ and $t\mapsto\gamma(c\,t)$, where $c>0$, have by definition the same end.

To make it more concrete: if we choose a point $\mathbf{T}\in\mathit{Pos}^{n+1}_1$, we get all the ideal points as the ends of the geodesics $\gamma$ such that $\gamma(0)=\mathbf{T}$, and two such geodesics $\gamma_1, \gamma_2$ have the same end iff $\gamma_1'(0)=c\,\gamma_2'(0)$ for some $c>0$.

Our ``data-space'' $\mathbb{RP}^n$ can be seen as a part of this ideal boundary.
Namely, a geodesic \eqref{pos-geod} ends at $[\mathbf{\tilde x}]\in\mathbb{RP}^n$ iff $\mathbf{\tilde x}$ is an eigenvector of the traceless matrix $\mathbf{T}^{-1}\mathbf{V}$ with a negative eigenvalue and if $\mathbf{T}^{-1}\mathbf{V}$ has only one more eigenvalue, with multiplicity $n$ (the other eigenspace is then the $\mathbf{T}$-orthogonal complement of the line $[\mathbf{\tilde x}]$).

To put it more explicitly, the geodesic $\gamma$ given by
$$\gamma(0)=\mathbf{T},\quad\gamma'(0)=-\on{grad}_\mathbf{T}\log\bigl(\mathbf{\tilde x}^T\mathbf{T}\,\mathbf{\tilde x}\bigr)
=
\frac{1}{n+1}\,\mathbf{T} - \mathbf{T\, \frac{\tilde x\tilde x^\mathit{T}}{\tilde x^\mathit{T}T\tilde x}\,T} 
$$
ends at $[\mathbf{\tilde x}]\in\mathbb{RP}^n$ (for this $\gamma$ we have $\Vert\gamma'(0)\Vert = \sqrt{n/(n+1)}$).

\subsubsection{Maximum likelihood, equilibrium of unit forces, and geodesic convexity}\label{sec:mvcauchy:end}

For a point $\mathbf{x}\in\mathbb{R}^n$ the negative log likelihood is, up to an inessential factor $(n+1)/2$, 
$$\log\bigl(\mathbf{\tilde x}^T\mathbf{T}\,\mathbf{\tilde x}\bigr).$$
This is (again up to an inessential factor $\sqrt{n/(n+1)}$) the ``regularized distance between $\mathbf{T}$ and the point at infinity $[\mathbf{\tilde x}]\in\mathbb{RP}^n\,$'', which is more properly called the Busemann function. It is defined by
$$B_{[\mathbf{\tilde x}]}(\mathbf{T})=
\lim_{t\to\infty}d\bigl(\mathbf{T},\gamma(t)\bigr) - d\bigl(\mathbf{T}_0,\gamma(t)\bigr)$$
where $\gamma$ is any geodesic ending in $[\mathbf{\tilde x}]$ and $\mathbf{T}_0\in\mathit{Pos}^{n+1}_1$ is arbitrarily chosen; it is well-defined only up to an additive constant, because of the dependence on $\mathbf{T}_0$.

The gradient of $B_{[\mathbf{\tilde x}]}(\mathbf{T})$ is the unit vector field tangent to the geodesics starting at $[\mathbf{\tilde x}]$. This gives us the ``equilibrium of unit forces pointing at the data-points'' criterion for the maximum likelihood parameter $\mathbf{T}$. Moreover, the Busemann function is geodesically convex, due to the fact that $\mathit{Pos}^{n+1}_1$ is non-positively curved and complete.
 We can therefore use the gradient descent along geodesics to find the maximum likelihood parameter.
 
\subsubsection{The covariant Hessian of the log likelihood and a safe step of the geodesic gradient descent}\label{sec:mvcauchy-hessian}
We already know that if $\gamma$ is a geodesic and if $\mathbf{\tilde x}\in\R^{n+1}\setminus\{0\}$, then the function
$$h(t):=\log\bigl(\mathbf{\tilde x}^T\gamma(t)\,\mathbf{\tilde x}\bigr)$$
 is convex. Let us now find and bound its 2nd derivative. Using the action of $\mathit{SL}(n+1)$ we can suppose that $\gamma(0)=\mathbf{I}$.
We then have, for $\mathbf{V}=\gamma'(0)$,
\begin{align*}
h(t)&=
\log\bigl(\mathbf{\tilde x}^T\exp(t\,\mathbf{V})\,\mathbf{\tilde x}\bigr)\\
&=\log(\mathbf{\tilde x}^T\mathbf{\tilde x})+
t\, \on{Tr}\mathbf{PV}
+ \tfrac12t^2\,\on{Tr}\bigl(\mathbf{PV}^2 - \mathbf{PVPV}\bigr)
 + O(t^3)
\end{align*}
where
$$\mathbf{P=\frac{\;\;\tilde x\,\tilde x^\mathit{T}}{\tilde x^\mathit{T}\tilde x}}$$
is the orthogonal projector onto the line $[\mathbf{\tilde x}]\subset\R^{n+1}$.

Since
\begin{equation}\label{mvcauchy-hess-bound}
0\leq\on{Tr}\bigl(\mathbf{PV}^2 - \mathbf{PVPV}\bigr)\leq\on{Tr}\mathbf{V}^2=\Vert\mathbf{V}\Vert^2
\end{equation}
(and since we can shift $t$ by any constant), we have for all $t\in\R$
$$0\leq h''(t)\leq\Vert\gamma'(t)\Vert^2 = \Vert\gamma'(0)\Vert^2.$$

If we have datapoints $\mathbf{\tilde x}_i$, $i=1,\dots, N$,  we would like to minimize the negative log likelihood (divided by $N$ to get nicer formulas)
\begin{equation}\label{mvcauchy-llik}
\ell(\mathbf{T})=\frac1N\sum_{i=1}^N\log\bigl(\mathbf{\tilde x}_i^T\mathbf{T}\,\mathbf{\tilde x}_i\bigr)
\end{equation}
via geodesic gradient descent. If $\gamma$ is a geodesic and if
$$g(t)=\ell(\gamma(t)),$$
we  have  $g''(t)\leq \Vert\gamma'(0)\Vert^2$. Moreover, if $\gamma'(0)$ is minus the gradient of $\ell$ at $\gamma(0)$, then $g'(0)=-\Vert\gamma'(0)\Vert^2$; therefore, a safe step for the descent is $1$.
 In other words, if we are at a point $\mathbf{T}_k$ and if $\mathbf{V}_k$ is the gradient%
\footnote{For completeness: the gradient of $\ell$ is
$$- \frac{1}{n+1}\mathbf{T} + \frac1N\mathbf{T}\,\sum_{i=1}^N\mathbf{\frac{\tilde x_{\mathit{i}}\tilde x_{\mathit{i}}^\mathit{T}}{\tilde x_{\mathit{i}}^\mathit{T}T\tilde x_{\mathit{i}}}}\,\mathbf{T} $$}
of $\ell$ at $\mathbf{T}_k$, we use the geodesic $\gamma$ given by $\gamma(0)=\mathbf{T}_k$, $\gamma'(0)=-\mathbf{V}_k$ (see \eqref{pos-geod}),  move to the point
$$\mathbf{T}_{k+1} = \gamma(1)$$
and iterate. At each step we are guaranteed to decrease $\ell$ at least by $\frac1{2}\Vert\mathbf{V}_k\Vert^2$.

In practice, the step size $1$, coming from the upper bound in \eqref{mvcauchy-hess-bound}, is somewhat too conservative. Namely, the Hessian 
$$\mathbf{V}\mapsto\on{Tr}\bigl(\mathbf{PV}^2 - \mathbf{PVPV}\bigr)$$
is of the form $(\mathbf{V}, \mathcal{P}\mathbf{V})$, where $\mathcal{P}$ is a rank-$n$ orthogonal projector in the space of symmetric traceless matrices. Since the space of symmetric traceless $(n+1)\times(n+1)$ matrices has dimension $(n^2+3n)/2$, it suggests the improved step size
$$\frac{n+3}{2}$$
(since $\mathbf{P}$'s would tend to average to a multiple of $\mathbf{I}$ and that multiple can be found by computing the trace).
I haven't found a proof that this step size is safe, but in practice it works well.

\subsubsection{Existence and uniqueness}\label{sec:mvcauchy-exun}
 For completeness, let us use convexity to reprove that, under general conditions, the maximum likelihood problem has a unique solution (a fact known since \cite[Theorem 4.1]{KT}). More precisely: 
\medskip

\emph{If we have points $\mathbf{\tilde x}_i\in\R^{n+1}\setminus\{0\}$, $i=1,\dots, N$, such that $N\geq n+2$ and such that no $n+1$ of the $\mathbf{\tilde x}_i$'s are linearly dependent, then the maximum likelihood problem
$$\mathbf{\hat T}=\argmin_{\mathbf{T}\in\mathit{Pos}^{n+1}_1}
\sum_{i=1}^N\log\bigl(\mathbf{\tilde x}_i^T\mathbf{T}\,\mathbf{\tilde x}_i\bigr)$$
has a unique solution.}%
\footnote{If $N<n+2$ then this cannot be true: the subgroup of $\mathit{SL}(n+1)$ preserving the lines $[\mathbf{\tilde x}_i]\subset\R^{n+1}$ is then closed and non-compact, and at the same time it would have to preserve $\mathbf{\hat T}$, i.e.\ be in a conjugate of $\mathit{SO}(n+1)$, giving a contradiction.

The hypothesis can be weakened to allow some linear dependencies, but then the statement becomes somewhat combinatoric and we shall omit it (see \cite{KT}).}
\begin{proof}
It is sufficient to prove that for any geodesic $\gamma$ the function
$$h(t)= \sum_{i=1}^N\log\bigl(\mathbf{\tilde x}_i^T\gamma(t)\,\mathbf{\tilde x}_i\bigr)$$
is strictly convex and that
$$\lim_{t\to\infty}h(t)=\lim_{t\to-\infty}h(t)=\infty.$$
In fact, it is enough to show just that $\lim_{t\to\infty}h(t)=\infty$ for any $\gamma$: then also $\lim_{t\to-\infty}h(t)=\infty$ by reversing the direction of $\gamma$, and since $h$ is real analytic, convex, and $\lim_{t\to\pm\infty}h(t)=\infty$, it then must be strictly convex.

By acting with a suitable element of $\mathit{SL}(n+1)$ we can suppose that 
$$\gamma(0)=\mathbf{I},\quad\gamma'(0)=\mathbf{D},$$
where $\mathbf{D}$ is a diagonal traceless matrix
$$\mathbf{D}=\on{diag}(\lambda_1,\dots,\lambda_{n+1}),\quad
\lambda_1\leq\dots\leq\lambda_{n+1},\ \sum_j\lambda_j=0,\ \sum_j\lambda_j^2>0.$$
This gives us
$\gamma(t)=\exp(t\,\mathbf{D})$, and thus
$$h(t)=\sum_{i=1}^N\log\Bigl(\sum_{j=1}^{n+1}c_{ij}^2 \exp(\lambda_j t)\Bigr)$$
where $c_{ij}$ are the components of the vector $\mathbf{\tilde x}_i$, and so we have the asymptotic
$$h(t)\sim \Bigl(\sum_{i=1}^N\lambda_{j(i)}\Bigr)\, t\qquad(\text{for }t\to\infty)$$
where
$$j(i)=\max\{j \mid c_{ij}\neq 0\}.$$
The hypothesis on $\mathbf{\hat x}_i$'s implies
$$\Bigl(\sum_{i=1}^N\lambda_{j(i)}\Bigr)\geq \Bigl(\sum_{j=1}^{n+1}\lambda_j\Bigr) +\lambda_{n+1} = \lambda_{n+1} >0,$$
therefore $\lim_{t\to\infty}h(t)=\infty$, as we wanted to prove.
\end{proof}

\begin{samepage}  
\subsection{Cauchy-like distributions on matrices}
\subsubsection{All $n\times m$ matrices}
Let us start with the distribution
\end{samepage}
$$f_{m,n}(\mathbf{X})\propto\det(\mathbf{I}_m + \mathbf{X}^T \mathbf{X})^{-(m+n)/2} $$
on the space $\mathit{Mat}(n,m)$ of $n\times m$ matrices, then we apply to it all the translations and the $GL(n)\times GL(m)$ action by left and right multiplication. We get the family of distributions
$$f_{m,n}(\mathbf{X}\mid\mathbf{S}_1,\mathbf{S}_2,\mathbf{B})\propto\det\bigl(\mathbf{S}_1 + (\mathbf{X}-\mathbf{B})^T \mathbf{S}_2^{-1}(\mathbf{X}-\mathbf{B})\bigr)^{-(m+n)/2}, $$
where the parameters $\mathbf{S}_1,\mathbf{S}_2$ are positive definite symmetric matrices and $\mathbf{B}\in\mathit{Mat}(n,m)$. This is the matrix-variate $t$-distribution with 1 degree of freedom \cite{GN}.

Again it is convenient to parametrize this family a bit differently. Let $\mathbf{\tilde X}$ be the $(m+n)\times m$-matrix obtained from $\mathbf{X}$ by appending the unit matrix $\mathbf{I}_m$ at the bottom. Then the same family of distributions can be parametrized by $\mathbf{T}\in\mathit{Pos}_{m+n}^1$:
\begin{equation}\label{matr-cauchy-inv}
f_{m,n}(\mathbf{X}\mid\mathbf{T})\propto\det(\mathbf{\tilde X}^T\mathbf{T}\,\mathbf{\tilde X})^{-(m+n)/2}.
\end{equation}
The normalization constant in \eqref{matr-cauchy-inv} is \emph{independent of the parameter} $\mathbf{T}$.

We extend $\mathit{Mat}(n,m)$ to $\mathit{Gr}(m+n, m)$ (the space of $m$-dim vector subspaces in $\R^{m+n}$) by attaching a set of measure 0: a matrix $\mathbf{X}\in\mathit{Mat}(n,m)$ corresponds to the graph (a subspace of $\R^{m+n}$) of the linear map $\mathbf{X}\colon\R^m\to\R^n$. Our family of distributions, living on $\mathit{Gr}(m+n, m)$, is closed under the action of $\mathit{SL}(m+n)$ on $\mathit{Gr}(m+n, m)$: the action of $\mathbf{A}\in \mathit{SL}(m+n)$ replaces the parameter $\mathbf{T}$ by $\mathbf{A}^T\mathbf{T}\,\mathbf{A}$.

Everything now works exactly as for the multi-variate Cauchy distribution (which is a special case, with $m=1$): the distributions are parametrized by $\mathit{Pos}_{m+n}^1$ and a geodesic $\gamma$ in $\mathit{Pos}_{m+n}^1$, as given by \eqref{pos-geod}, ends at $W\in\mathit{Gr}(m+n, m)$ iff the traceless matrix $\mathbf{T}^{-1}\mathbf{V}$ has only 2 eigenvalues and $W$ is the eigenspace corresponding to the negative one (the other eigenspace is then the $\mathbf{T}$-orthogonal complement of $W$). 

The safe step $1$ for the geodesic gradient descent still works. The projector $\mathcal P$ (\S\ref{sec:mvcauchy-hessian}) now projects onto the space of symmetric matrices which, when seen as quadratic forms, vanish both on $W$ and on $W^\perp$. Is has rank $mn$, which suggests the improved step size
$$\frac{(m+n)(m+n+1)-2}{2mn}.$$

\subsubsection{Symmetric and skew-symmetric matrices}
There are similar families of  probability densities on the spaces of symmetric and skew-symmetric $n\times n$ matrices:
$$f^\text{sym}_{n,n}(\mathbf{X})\propto\det\Bigl(\mathbf{I}_n+ \bigl(\mathbf{S}^{-1}(\mathbf{X}-\mathbf{B})\bigr)^2\Bigr)^{-(n+1)/2}$$
$$f^\text{skew}_{n,n}(\mathbf{X})\propto\det\Bigl(\mathbf{I}_n- \bigl(\mathbf{S}^{-1}(\mathbf{X}-\mathbf{B})\bigr)^2\Bigr)^{-(n-1)/2}.$$
The parameters are a symmetric positive definite matrix $\mathbf{S}$ (in both cases) and a symmetric or skew-symmetric $\mathbf{B}$ respectively.

We can still rewrite these distributions as
$$f^\text{sym}_{n,n}(\mathbf{X})\propto\det(\mathbf{\tilde X}^T\mathbf{T}\,\mathbf{\tilde X})^{-(n+1)/2}$$
$$f^\text{skew}_{n,n}(\mathbf{X})\propto\det(\mathbf{\tilde X}^T\mathbf{T}\,\mathbf{\tilde X})^{-(n-1)/2},$$
with the normalization constant independent of $\mathbf{T}$. However, now the parameter $\mathbf{T}\in\mathit{Pos}_{2n}^1$ needs to satisfy an extra condition. If we set
$$\mathbf{H}^\text{skew}=
\begin{bmatrix}
\mathbf{0}_n & \mathbf{I}_n\\
-\mathbf{I}_n & \mathbf{0}_n 
\end{bmatrix},
\qquad
\mathbf{H}^\text{sym}=
\begin{bmatrix}
\mathbf{0}_n & \mathbf{I}_n\\
\mathbf{I}_n & \mathbf{0}_n 
\end{bmatrix},
\qquad
$$
then $\mathbf{T}$ needs to satisfy 
$$(\mathbf{H}^\text{skew}\mathbf{T})^2 = -\mathbf{I}_{2n}$$
 in the \emph{symmetric} case, and 
$$(\mathbf{H}^\text{sym}\mathbf{T})^2 = \mathbf{I}_{2n}$$
in the \emph{skew-symmetric} case.

Nonetheless, these conditions give
 totally geodesic submanifolds of $\mathit{Pos}_{2n}^1$, so from the calculation point of view nothing really changes.
 
In place of the full Grassmannian $\mathit{Gr}(2n,n)$ we now get the (symplectic or symmetric, in the symmetric and skew-symmetric case respectively) Lagrangian Grassmannian -- the space of $n$-dim subspaces of $\mathbb R^{2n}$ on which the bilinear form given by $\mathbf{H}^\text{skew}$ (or $\mathbf{H}^\text{sym}$) vanishes.%
\footnote{To be precise, in the skew-symmetric case we are getting only one of the two connected components of the symmetric Lagrangian Grassmannian.} The group $SL(2n)$ gets replaced by the subgroup preserving either $\mathbf{H}^\text{skew}$ or $\mathbf{H}^\text{sym}$, i.e.\ by $\mathit{Sp}(2n)$ or by $\mathit{SO}(n,n)$.

\subsection{Probability distributions parametrized by symmetric spaces}\label{sec:general}
Let us now explain the general story.

Let $G$ be a non-compact semisimple Lie group and let $P\subset G$ be a parabolic subgroup. Let $K\subset G$ be a maximal connected compact subgroup ($K$ is unique up to conjugation). Then $K$ acts transitively on $G/P$, so we have a unique $K$-invariant probability density $\mu_K$ on $G/P$. If we transport $\mu_K$ by the action of an element $g\in G$, we get the $gKg^{-1}$-invariant density $\mu_{gKg^{-1}}$.

These densities on $G/P$ are parametrized by $G/K$. The manifold $G/K$ has a natural $G$-invariant Riemannian metric (it is a symmetric space of non-compact type), coming from the Killing form on the Lie algebra $\mathfrak{g}$ restricted to $\mathfrak{k}^\perp\subset\mathfrak{g}$. If $G$ is simple then this metric is, up to a constant multiple, the Fisher information metric of our family of probability densities (since $\mathfrak{k}^\perp$ is an irreducible $\mathfrak{k}$-module and so, up to a multiple, it admits a unique $\mathfrak{k}$-invariant inner product; for non-simple $G$ we might need to use different multiples for different factors of $G$).

One can then check that there is an embedding of $G/P$ to the ideal boundary of $G/K$ such that the negative log likelihoods are, up to a constant multiple, the ``regularized distances'' (i.e.\ Busemann functions) to the datapoints (which are points in $G/P$). We thus have the ``equilibrium of unit forces'' picture for the maximum likelihood point of $G/K$. Moreover, since $G/K$ has non-positive curvature and is complete and simply connected, the Busemann functions are geodesically convex.

For applications we might like to identify $G/P$, up to a subset of measure 0, with a vector space (say $\R^n$), and demand our family of densities to be closed under translations in the vector space. This can be achieved as follows.
Let $N\subset G$ be the nilpotent radical of $P$. Then $N$ has a ``big'' orbit $O\subset G/P$: the action of $N$ on $O$ is free and the complement of $O$ is of measure 0. We are thus done if $N$ is abelian, as in that case $N$ a vector space and we can identify $O$ with $N$ by choosing a point in $O$.

The classification of \emph{complex} pairs $P\subset G$ (with simple $G$) with this property is done in \cite{RRS}:
\begin{enumerate}[label={(\alph*)}]
\item $G=\mathit{SL}(n)$, $G/P$ is a Grassmannian $\mathit{Gr}(n, k)$,
\item $G=SO(n)$, $G/P\subset\mathbb{CP}^{n-1}$ is the corresponding projective quadric,
\item $G=SO(2n)$ or $Sp(2n)$, $G/P$ is the corresponding Lagrangian Grassmannian,
\item two more exceptional cases, one with $G=E_6$ and one with $G=E_7$.
\end{enumerate}
In our examples we simply took suitable real forms of the first 3 options: (a) corresponds to matrices $(n-k)\times k$, (b)\ to the ``conformal'' distributions on $\R^{n-2}$, and (c)\ to skew-symmetric and symmetric $n\times n$ matrices respectively. We didn't take all possible real forms (e.g\ in case (b) any $SO(n-k,k)$, with both $k$ and $n-k$ nonzero, admits a real form of $P$, but we considered only $k=1$; in all the cases we could also take the complex $P\subset G$ and see them as real Lie groups), but the examples we considered give a general idea about the probability densities that arise in this way.

\newpage
\appendix

\section{Numerical experiments and practical issues}\label{sec:appendix}
Let us now give a few numerical examples. There is certainly nothing new about the idea of using the maximum likelihood estimation (MLE) with the multivariate Cauchy distribution as a robust estimator of location and of scatter. The main interest is in seeing how our geodesic gradient descent (GGD) method performs. Let us start with a few general remarks.

\begin{itemize}
\item The gradient of the loss function $\ell$ (see \eqref{mvcauchy-llik}) has length smaller than 1 everywhere. It might therefore seem  that the GGD could take a lot of time if the starting point is not well chosen. Fortunately, ``nothing is far in the space $\mathit{Pos}_{n+1}^1$'': the (geodesic) distances grow only logarithmically in the entries of the matrices (cf.\ \eqref{pos-geod}). We always start the descent at the identity matrix. Nonetheless, \emph{numerical} stability might be an issue.
\item While the GGD typically converges rapidly, for some datasets it's not the case. This simply means (because of convexity) that the argmin is ill-conditioned (i.e.\ that there are points rather far from the argmin where the value of $\ell$ is close to the minimum) and so it's unstable under small changes in the dataset. In other words, in this case the MLE makes no statistical sense. Hence the descent algorithm runs always only for a short time, and either gives the estimate or a warning.
\item
In the examples we will take iid samples $\mathbf{x}_i$ ($i=1,\dots,N$) from a probability measure $P$ on $\R^n$. The function
$$\ell_P(\mathbf{T})=\on{E}_P\bigl[\log(\mathbf{\tilde x^\mathit{T}T\tilde x})-\log(\mathbf{\tilde x^\mathit{T}\tilde x})\bigr]$$
(where $\mathbf{x}$ is the $P$-distributed random variable; the subtracted term just ensures existence of the expectation) is geodesically convex. In good cases (which are generic -- see \S\ref{sec:mvcauchy-exun}) the function $\ell_P$ goes to $\infty$ at the ideal boundary and is strictly geodesically convex, hence
$$\mathbf{T}_P=\argmin\ell_P$$
is well-defined.
Our sequence of estimators
$$\mathbf{\hat T}_N = \argmin_{\mathbf{T}\in\mathit{Pos}^{n+1}_1}
\frac1N\sum_{i=1}^N\log\bigl(\mathbf{\tilde x}_i^T\mathbf{T}\,\mathbf{\tilde x}_i\bigr)$$
then converges to $\mathbf{T}_P$ and is asymptotically normally distributed \cite{Huber}. The covariance matrix of this normal distribution%
\footnote{This covariance matrix is \cite[Corollary 6.7]{Huber}
$$\tfrac1N\,\bigl(\on{Hess}_{\mathbf{T}_P}\ell_P\bigr)^{-1} \on{E}_P\bigl[\bigl(\on{grad}_{\mathbf{T}_P}\log(\mathbf{\tilde x^\mathit{T}T\tilde x})\bigr)^{\otimes2}\bigr] 
\bigl(\on{Hess}_{\mathbf{T}_P}\ell_P\bigr)^{-1}.$$
If we replace $P$ by the empirical distribution given by a dataset (in particular, replace $\ell_P$ by $\ell$), it gives us a first estimate about the precision of $\mathbf{\hat T}_N$. We do not do this analysis in what follows.
}
 is large when the Hessian of $\ell_P$ at $\mathbf{T}_P$ is close to being degenerate (and the asymptotic regime then kicks in only for very large $N$'s). Hence, as we see again, using the MLE is not a good idea it this case.
\end{itemize}

An implementation of the algorithm can be found \href{https://gist.github.com/mrwth/0e7a3885d970d9d87b0123afad2227d0}{\texttt{here}}.

\subsection{1d examples}
From each distribution that we consider we take an iid sample of size 1000 and run the GGD algorithm. We repeat it 1000 times, each time with a new sample.

Let us start with the standard normal distribution $N(0,1)$. We get the following results.
$$\includegraphics[scale=0.8]{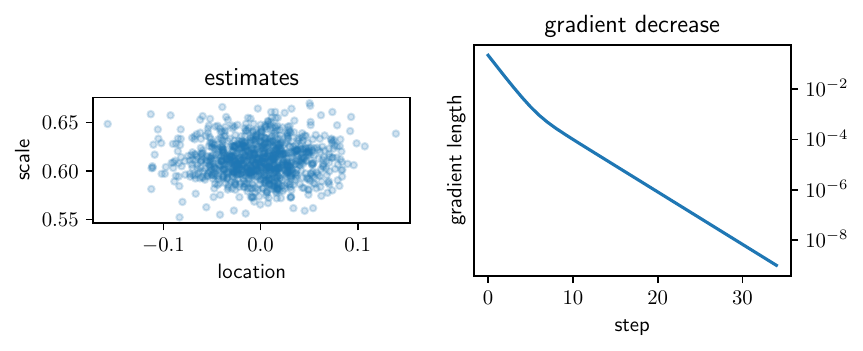}$$
The ``location'' and the ``scale'' parameters correspond to $u$ and $v$ in \eqref{1dcauchy}. They are reasonably close to the limit values $0$ and $0.612$. The graph on the right shows how the gradient was behaving during the descent (shown only for the last out of 1000 runs).

Now the same, with a mixed distribution: 90\%  $N(0,1)$ and 10\%  $N(100, 100^2)$.
$$\includegraphics[scale=0.8]{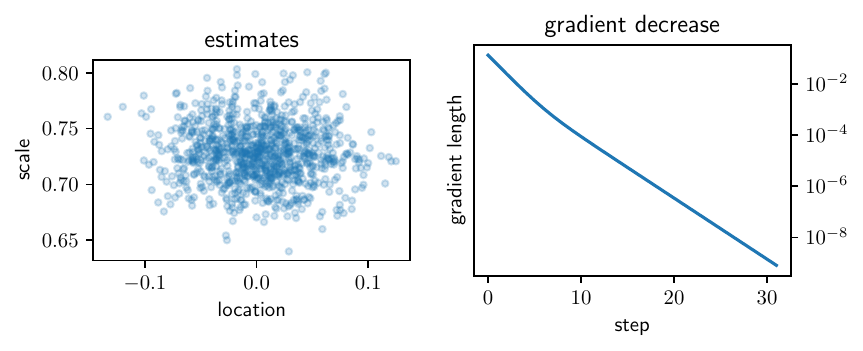}$$
As we see, the outliers didn't move the estimated location at all and slightly increased the estimated scale.

Two more distributions where the method works well: Cauchy (of course!) with $u=1000$ and $v=10$ 
$$\includegraphics[scale=0.8]{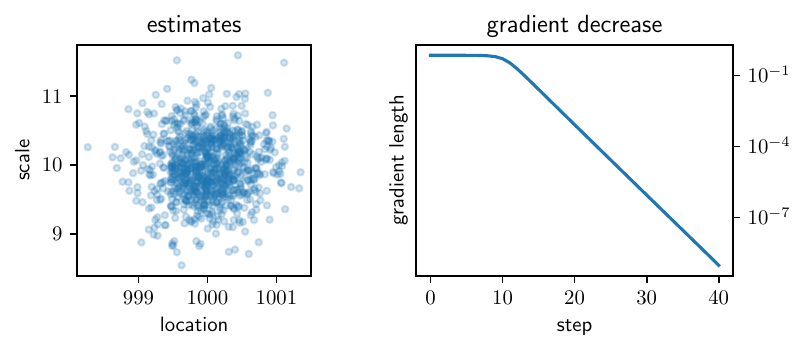}$$
(the graph on the right shows that it took a few steps until the gradient started decreasing fast -- that's because the starting point of the GGD is $u=0$, $v=1$)
and normal $N(10, 3^2)$ with 10\% contamination by Cauchy with $u=0$, $v=1$
$$\includegraphics[scale=0.8]{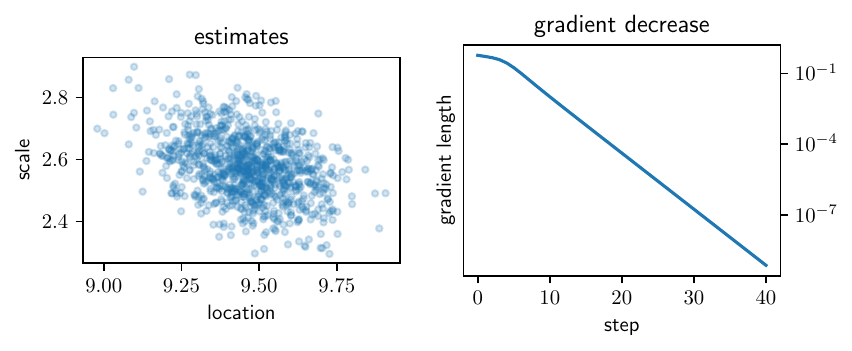}$$
(though here the contamination moved the location estimates away from 10 a bit).

It is easy to understand the cases where the 1d Cauchy MLE struggles. We have existence and unicity of $\mathbf{T}_P$ for all the probability measures $P$ on $\R\cup\{\infty\}=\mathbb{RP}^1$ except for those that have an atom of measure at least $1/2$. If we have 2 atoms with probability 1/2 then $\mathbf{T}_P$ (defined via argmin) is not unique -- they form an entire geodesic. In the other cases $\mathbf{T}_P$ doesn't exist, but if we extend the class of Cauchy distributions to include also $\delta$-functions, then it does exist and it is unique (it is $\delta$ at the position of the atom). Hence the MLE, and also our GGD, will struggle for $P$'s which are ``close'' to those bad 2-atomic measures, i.e.\ when we have 2 well-separated lumps with roughly the same probability (one of the lumps might be around $\infty$).

Here is what happens when we take a 50/50 mixture of $N(0, 10^2)$ and $N(300, 1)$ (with a geodesic showing up on the left pane):
$$\includegraphics[scale=0.8]{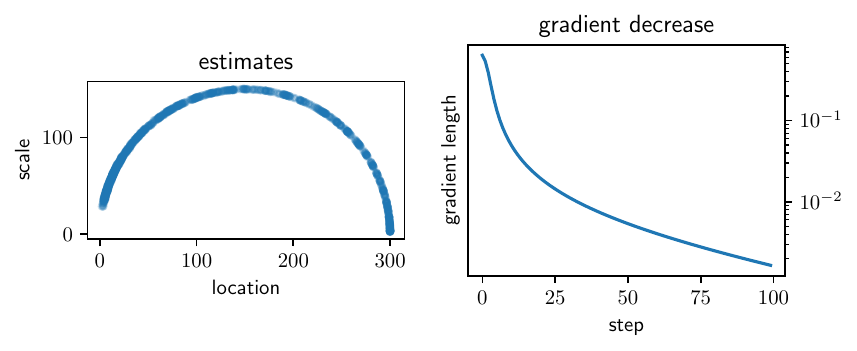}$$
Already the graph on the right (requiring just a single run) shows us that the MLE is not stable. For the mixture ratio  60/40, the method is starting to recover:
$$\includegraphics[scale=0.8]{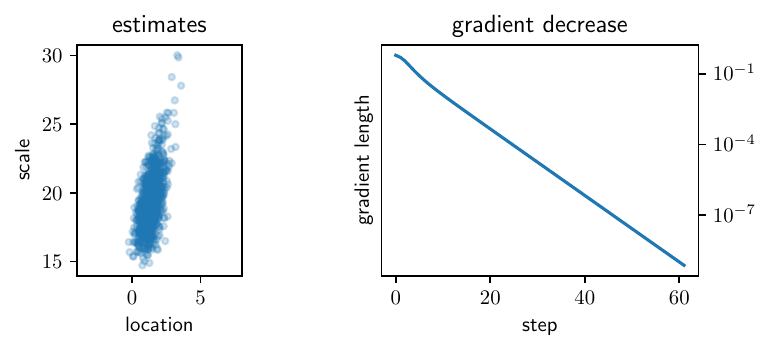}$$
Figuratively speaking, the Cauchy MLE is trying hard to decide which part of the data is formed by outliers, and it has hard time doing it in these cases.

\subsection{2d examples}
We shall still use 1000 datapoints and 1000 runs. We are now estimating 5 parameters, so we may expect a bit more noise.

Let us start with the normal distribution
$N\biggl(
\begin{bmatrix}
2\\-3
\end{bmatrix},
\begin{bmatrix}
2&1\\
1&3
\end{bmatrix}
\biggr)
$.
$$\includegraphics[scale=0.8]{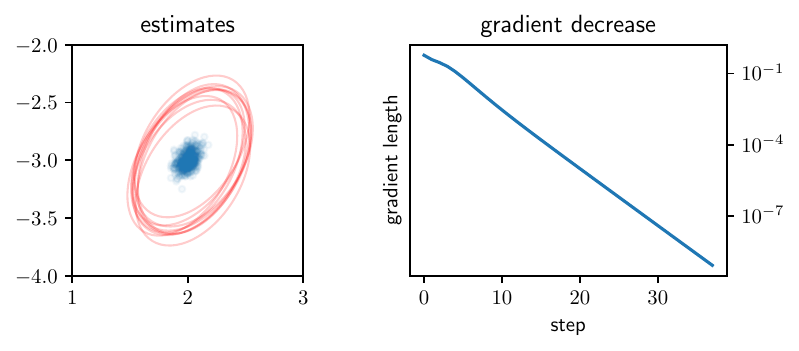}$$
On the left we show the estimated location $\mathbf{\hat b}$ for all the 1000 runs, and for the first  10  of them (to not clutter the picture too much) we show also the corresponding scatter $\mathbf{\hat S}$.

Let us add to it 10\% contamination by 
$N\biggl(
\begin{bmatrix}
20\\70
\end{bmatrix},
\begin{bmatrix}
15&8\\
8&9
\end{bmatrix}
\biggr)
$.
$$\includegraphics[scale=0.8]{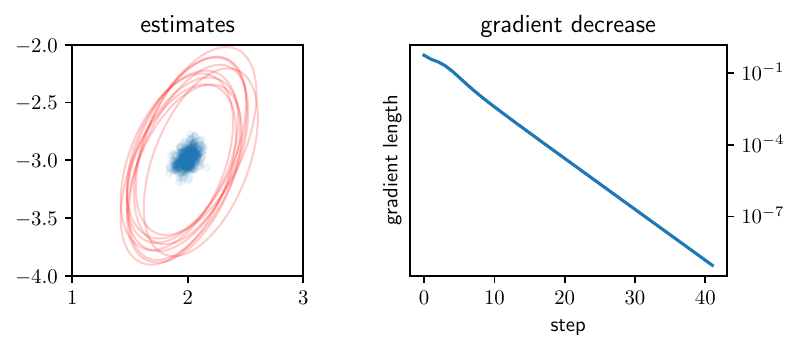}$$
As we can see, the location estimate doesn't change and the scatter estimate is slightly modified. The algorithm still runs fast.

\subsection{4d examples}

Let us start with $N(\bm{\mu},\bm{\Sigma})$
where
$$
\bm{\mu}=
\begin{bmatrix}
1\\
2\\
3\\
4
\end{bmatrix}
\qquad
\bm{\Sigma}=
\begin{bmatrix}
1&1&1&1\\
1&2&2&2\\
1&2&3&3\\
1&2&3&4
\end{bmatrix}
$$
We are estimating 14 coefficients. This time we shall do only one run, with a large sample of size $10^7$ (and not do a ``visual statistical analysis''). At this point we just want to illustrate the efficiency of the algorithm. 

The algorithm stopped after 34 steps (when the gradient size dropped below $10^{-9}$) and produced the estimate
$$
\mathbf{\hat b}=
\begin{bmatrix}
1.000\\
2.001\\
2.000\\
4.000
\end{bmatrix}
\qquad
\mathbf{\hat S}\propto
\begin{bmatrix}
1.000&1.000&1.000&0.999\\
1.000&2.000&2.000&2.000\\
1.000&2.000&3.000&3.000\\
0.999&2.000&3.000&3.999
\end{bmatrix}
$$
(the proportionality constant is chosen to get a matrix with the same determinant as $\bm{\Sigma}$).
After a contamination by 5\% of $N(\bm{\mu}',\bm{\Sigma}')$, where
$$
\bm{\mu}'=
\begin{bmatrix}
100&0&-100&0
\end{bmatrix}^T
\qquad
\bm{\Sigma}'=
500\,\bm{I}_4,
$$
the algorithm again took just 35 steps and produced
$$
\mathbf{\hat b}=
\begin{bmatrix}
1.001\\
2.000\\
2.999\\
3.998
\end{bmatrix}
\qquad
\mathbf{\hat S}\propto
\begin{bmatrix}
1.011& 0.910& 0.815& 0.908\\
0.910& 1.825& 1.823& 1.822\\
0.815& 1.823& 2.878& 2.738\\
0.908& 1.822& 2.738& 3.646
\end{bmatrix}
$$

\end{document}